\documentclass[12pt]{amsart}
\usepackage{amsmath,amsthm,latexsym,amscd,amsbsy,amssymb,amsfonts,fleqn,leqno,
euscript, pb-diagram}

\numberwithin{equation}{section}
\newtheorem{thm}{Theorem}[section]

\newtheorem{conj}[thm]{Conjecture}
\newtheorem{cor}[thm]{Corollary}
\newtheorem{qu}[thm]{Question}

\def\leukfrac#1/#2{\leavevmode
               \kern.1em
                \raise.9ex\hbox{\the\scriptfont0 ${}_#1$}
                \hskip -1pt\kern-.1em
                /\kern-.15em\lower.10ex\hbox{\the\scriptfont0 ${}_#2$}}

\theoremstyle{definition}

\theoremstyle{remark}

\theoremstyle{definition}

\newcommand{\bd}{\mathrm{bd}}

\begin{document}
%%%%%%%%%%% Begin Topmatter %%%%%%%%%%%%%%%%%

\title[Homogeneous metric $ANR$-compacta]
{Homogeneous metric $ANR$-compacta}

\author{V. Valov}
\address{Department of Computer Science and Mathematics,
Nipissing University, 100 College Drive, P.O. Box 5002, North Bay,
ON, P1B 8L7, Canada} \email{veskov@nipissingu.ca}

\date{\today}
\thanks{The author was partially supported by NSERC
Grant 261914-19.}


\begin{abstract}
This is a survey of most important results and unsolved problems about homogeneous finite-dimensional metric $ANR$-compacta. We also discuss some partial results and possible ways of solutions.
\end{abstract}
\maketitle\markboth{}{Homogeneous $ANR$,s}
%%%%%%%%%% End topmatter %%%%%%%%%%%%%%%%%%%%%

%%%%%%%%%%%%%%%%%%%%%%%%%%%%%%%%%%%%%%%%%%%%%%%%%%%%%%%%%%%%%%%%
%%%%%%%%%%%%%%%%%%%%%%%%%%%%%%%%%%%%%%%%%%%%%%%%%%%%%%%%%%%%%%%%

%%%%%%%%%%%%%%%%%% TABLE OF CONTENT %%%%%%%%%%%%%%%%%%%%%%%%%%%%%%%%%%

%\tableofcontents

%%%%%%%%%%%%%%%%%%%%%%%%%%%%%%%%%%%%%%%%%%%%%%%%%%%%%%%%%%%%%%%%%%%
%%%%%%%%%%%%%%%%%%%%%%%%%%%%%%%%%%%%%%%%%%%%%%%%%%%%%%%%%%%%%%%%%%%%%%

\section{Introduction}
%All spaces in this paper are assumed to be at least normal.

In this paper we will survey the most important results and unsolved problems concerning homogeneous finite-dimensional $ANR$-compacta, their relationship to each other, as well as possible ways of solutions. There are many interesting problems concerning this class of compacta. Definitely, the most important problem in this area is the Bing-Borsuk conjecture \cite{bb} stating that every compact homogeneous $ANR$-compactum of dimension $n$ is an $n$-manifold. This is true in dimensions 1 and 2, but still unknown for $n\geq 3$. Recall that an $n$-manifold is a separable metric space such that each point has a neighborhood homeomorphic to the Euclidean $n$-space $\mathbb R^n$.
It is well known that every compact connected $n$-manifold without boundary is homogeneous.
The importance of the Bing-Borsuk conjecture was confirmed by Jacobsche \cite{ja} that the $3$-dimensional Bing-Borsuk conjecture implies the celebrated Poincare conjecture proven by Perelman, see \cite{per}.
A special case of the Bing-Borsuk conjecture is the Busemann conjecture \cite{bu1}-\cite{bu2} stating that the Busemann $G$-spaces are manifolds. The Busemann conjecture has been proven for $G$-spaces of dimension $\leq 4$. For more information about the Busemann conjecture and related results, see \cite{hr}.

During the Spring Topology Conference in March 2018 J. Bryant and S. Ferry announced they constructed a counter-example to the Bing-Borsuk conjecture. This example is not published yet, so the Bing-Borsuk conjecture is still open.

Recall that a metric space $X$ is an absolute neighborhood retract (br., $ANR$) if for every embedding of $X$ as a closed subset of a metric space $M$ there exists a neighborhood $U$ of $X$ in $M$ and a retraction $r:U\to X$, i.e. a continuous map $r$ with $r(x)=x$ for all $x\in X$. A space $X$ said to be homogeneous if for any two points $x,y\in X$ there is a homeomorphism of $X$ onto itself taking $x$ to $y$.

Unless stated otherwise, all spaces are separable metric and all maps are continuous. Predominantly, reduced \v{C}ech homology $\check{H}_n(X;G)$ and cohomology $\check{H}^n(X;G)$ with coefficients from an abelian group $G$ will be considered. Singular homology and cohomology groups are denoted, respectively, by $H_k(X;G)$ and $H^k(X;G)$. By a dimension we mean the covering dimension $\dim$, the cohomological dimension with respect to a group $G$ is denoted by $\dim_G$.

\section{Homogeneous spaces and generalized Cantor manifolds}

Although the Bing-Borsuk conjecture is still open, it is known there are some common properties of homogeneous finite-dimensional $ANR$-spaces and Euclidean manifolds. Brouwer \cite{bro} proved a century ago that every $n$-manifold $X$ has the {\em invariance of the domain property}. This means that if $U,V$ are homeomorphic subsets of $X$, then $U$ is open if and only if $V$ is open. Cantor manifolds is another notion introduced by Urysohn \cite{ur} in 1925 as a generalization of Euclidean manifolds.
\begin{itemize}
\item An compact metric space $X$ is a {\em Cantor $n$-manifold} if $\dim X=n$ and $X$ can not be expressed as the union of two closed proper subsets whose intersection is of dimension $\leq n-2$.
\end{itemize}
One of the first results about homogeneous $ANR$-spaces is the following result of Lysko \cite{ly}:
\begin{thm}\cite{ly}
Every connected $n$-dimensional homogeneous $ANR$-compactum is a Cantor $n$-manifold and has the invariance of the domain property.
\end{thm}
 Lysko's result was generalized by Seidel \cite{se} for locally compact and locally homogeneous $ANR$'s. Recall that a space $X$ is {\em locally
homogeneous} provided for every two points $x,y\in X$ there exist neighbourhoods $U$ and $V$ of $x$ and $y$, respectively, and a
homeomorphism $h:U\to V$ such that $h(x)=y$.
 Krupski \cite{kr3} extended the first part of Lysko's result.
\begin{thm}\label{kr}\cite{kr3} Every homogeneous, locally compact and connected space of dimension $n$ is a Cantor $n$-manifold.
\end{thm}
 A new dimension, unifying both the covering and the cohomological dimension, was introduced in \cite{kktv}, and Theorem \ref{kr} was extended for this  dimension.
A sequence $\mathcal{K}=\{K_0,K_1,...\}$ of $CW$-complexes is called
a {\em stratum} for a dimension theory~\cite{dr} if
\begin{itemize}\item
$K_n\in AE(X)$, where $X$ is a metric space, implies both $K_{n+1}\in AE(X\times [0,1])$ and
$K_{n+j}\in AE(X)$ for all $j\geq 0$.
\end{itemize}
Here, $K_n\in AE(X)$ means that $K_n$ is an absolute extensor for
$X$. Given a stratum $\mathcal{K}$, the dimension function
$D_{\mathcal{K}}$ for a metrizable  space $X$ is defined as follows:
\begin{enumerate}
\item
$D_{\mathcal{K}}(X)=-1$ iff $X=\varnothing$;
\item $D_{\mathcal{K}}(X)\le n$ if
$K_n\in AE(X)$ for $n\ge 0$; if $D_{\mathcal{K}}(X)\le n$ and
$K_m\not\in AE(X)$  for all $m<n$, then $D_{\mathcal{K}}(X)= n$;
\item
 $D_{\mathcal{K}}(X)=\infty$ if $D_{\mathcal{K}}(X)\le n$ is not satisfied for any $n$.
\end{enumerate}

If $\mathcal{K}=\{\mathbb S^n\}_{n=0}^\infty$ is the sequence of all $n$-dimensional spheres, we obtain the covering dimension $\dim$.
If $G$ is a group and $\mathcal{K}=\{K(G,n)\}_{n=0}^\infty$, where $K(G,n)$ are the Eilenberg-MacLane complexes for $G$, then
the dimension $D_{\mathcal{K}}$ coincides with the cohomological dimension $\dim_G$.
We denote by $\mathcal D_{\mathcal K}^k$ the class of all metrizable spaces $X$ with
$D_{\mathcal{K}}(X)\leq k$.

Mazurkiewicz established another property of Euclidean spaces: any region in $\mathbb R^n$ cannot be cut by subsets of dimension $\leq n-2$
(a subset cuts if its compliment is not continuum-wise connected), see \cite{en}.
Inspired by this result, Hadjiivanov-Todorov \cite{ht} introduced the class of {\em Mazurkiewicz manifolds}.
This notion was generalized in \cite{kktv}.
\begin{itemize}
\item A normal space (not necessarily metrizable) $X$  is a {\em Mazurkiewicz manifold with respect to  $\mathcal{C}$}, where
$\mathcal{C}$ is a class of spaces,
if for every two closed, disjoint subsets $X_0,X_1\subset X$, both
having non-empty interiors in $X$, and every $F_\sigma$-subset
$F\subset X$ with $F\in\mathcal{C}$, there exists a continuum $K$ in
$X\setminus F$ joining $X_0$ and $X_1$. If in that definition $K$ can be chosen to be an arc, we say that $X$ is a {\em Mazurkiewicz arc manifold with respect to  $\mathcal{C}$} \cite{tv}.
\end{itemize}
Next theorem provides the strongest property of type Cantor manifolds possessed by homogeneous spaces.
\begin{thm}\cite{kv} Let $\mathcal{K}$ be a stratum and $X$ be a homogeneous locally compact, locally connected metric space.
Then every region $U\subset X$ with $D_{\mathcal{K}}(U)=n$ is a Mazurkiewicz manifold with respect to the class $\mathcal D_{\mathcal K}^{n-2}$.
\end{thm}

Alexandroff \cite{ps} introduced another property which is possessed by compact closed $n$-manifolds, to so-called continua $V^n$. Here is the general notion of  Alexandroff manifold, see \cite{tv}:
\begin{itemize}
\item A connected space $X$ is an
{\em Alexandroff manifold with respect to a given class $\mathcal{C}$} of spaces if for every two disjoint closed subsets $X_0,X_1$ of $X$, both having non-empty interiors, there exists an open cover $\omega$ of $X$ such that no partition $P$ between $X_0$ and $X_1$ admits an $\omega$-map onto a space $Y\in\mathcal{C}$. The Alexandroff {\em continua $V^n$} are Alexandroff manifolds with respect to the class of all spaces $Y$ with $\dim Y\leq n-2$.
\end{itemize}
Recall that a partition between two disjoint sets $X_0,X_1$ in $X$ is a closed set $F\subset X$ such that $X\setminus F$ is the union of two open disjoint sets $U_0,U_1$ in $X$ with $X_0\subset U_0$ and $X_1\subset U_1$. An $\omega$-map $f:P\to Y$ is such a map that $f^{-1}(\gamma)$  refines $\omega$ for some open cover $\gamma$ of $Y$. A cohomological version of $V^n$-continua was considered in \cite{s}, see also \cite{ktv} and \cite{t} for a subclass of the $V^n$-continua:
\begin{itemize}
\item A compactum $X$ is a {\em $V^n_G$-continuum} \cite{s}, where $G$ is a given group, if for every open disjoint subsets $U_1,U_2$ of $X$ there is an open cover
$\omega$ of $X_0=X\setminus (U_1\cup U_2)$ such that any partition $P$ in $X$ between $U_1$ and $U_2$ does not admit an $\omega$-map $g$ onto a space $Y$ with $g^*:\check{H}^{n-1}(Y;G)\to \check{H}^{n-1}(P;G)$ being a trivial homomorphism. If, in addition, there is also an element
$\gamma\in\check{H}^{n-1}(X_0;G)$ such that for any partition $P$ between $U_1$ and $U_2$ and any $\omega$-map $g$ of $P$ into a space $Y$ we have $0\neq i_P^*(\gamma)\in g^*(\check{H}^{n-1}(Y;G))$, where $i_P$ is the embedding $P\hookrightarrow X_0$, $X$ is called a {\em strong $V^n_G$-continuum} \cite{vv3}. A relative version of $V^n_G$-continua was considered in \cite{tv}.
\end{itemize}
Because $\check{H}^{n-1}(Y;G)=0$ for every group $G$ and every compact space $Y$ with $\dim Y\leq n-2$, all $V^n_G$-continua are $V^n$.

The following question is one of the remaining problems in that direction, see \cite{ktv}:
\begin{qu} Let $X$ be a homogeneous $ANR$-continuum and $G$ a group.
\begin{itemize}
\item[(1)] Is $X$ a $V^n$-continuum provided $\dim X=n$?
\item[(2)] Is $X$ a $V^n_G$-continuum provided $\dim_GX=n$?
\end{itemize}
\end{qu}

Karassev \cite{ka} provided a positive answer to Question 2.4(1) if $X$ is strongly locally homogeneous. Recall that a space is strongly locally homogeneous if every point $x\in X$ has a local basis of open sets $U$ such that for every $y,z\in U$ there is a homeomorphism $h$ on $X$ with $h(y)=z$ and $h$ is identity on $X\setminus U$. Every connected strongly locally homogeneous space is homogeneous.
Question 2.4 has also a positive answer if additionally $H^n(X;G)\neq 0$, see \cite[Corollary 1.2]{vv3}.

\section{Homology manifolds and the Modified Bing-Borsuk conjecture}

Topological $n$-manifolds $X$ have the following property: For every $x\in X$ the groups $H_k(X,X\setminus\{x\};\mathbb Z)$ are trivial if $k<n$   and $H_n(X,X\setminus\{x\};\mathbb Z)=\mathbb Z$. A space with this property is said to be a {\em $\mathbb Z$-homology $n$-manifold}. A {\em generalized $n$-manifold} is a locally compact $n$-dimensional $ANR$-space which is a $\mathbb Z$-homology $n$-manifold.
Every generalized $(n\leq 2)$-manifold is known to be a topological $n$-manifold \cite{w}. On the other hand, for every $n\geq 3$ there exists a generalized $n$-manifold $X$ such that $X$ is not locally Euclidean at any point, see for example \cite{can}. Any $n$-dimensional resolvable space is a generalized $n$-manifold.
\begin{itemize}
\item Recall that a space $X$ is {\em resolvable} if there exists a proper surjective map $f:M\to X$, where $M$ is a manifold, such that for each $x\in X$, $f^{-1}(x)$ is contractible in any neighborhood of itself in $M$. Such maps are said to be {\em cell-like maps}. The following property utilized by Cannon \cite{can1} plays a key role in geometric topology:
        A space $X$ has {\em the disjoint disks property} if arbitrary maps $f,g:\mathbb B^2\to X$ from the 2-dimensional disk $\mathbb B^2$ into $X$ can be approximated, arbitrary closely, by maps having disjoint images.
\end{itemize}
The characterization of manifolds culminated with the Edwards's theorem:
\begin{thm}\cite{ed}\label{main}
Topological $n$-manifolds, $n\geq 5$, are precisely the $n$-dimensional resolvable spaces having the disjoin disks property.
\end{thm}
Daverman-Repov\v{s} \cite{dr} provided a similar characterization of $3$-manifolds replacing the disjoint disks property by the so-called {\em simplicial spherical approximation property}, see also \cite{dr1}. Only partial results in dimension $4$ are known, see \cite{bdvw}, \cite{dr}.
Theorem \ref{main} shows the importance of the question as to whether generalized manifolds are resolvable, and this question had been for a long time, see \cite{can}.
\begin{conj}[Resolution conjecture] Every generalized $(n\geq 3)$-manifold is resolvable.
\end{conj}
In dimension $3$ the Resolution conjecture implies the Poincare conjecture \cite{re}, and only partial cases are known. In higher dimensions the Resolution Conjecture is false.  According to Bryant-Ferry-Mio-Weinberger \cite{bfmw1} there exist non-resolvable generalized $n$-manifolds for every $n\geq 6$. The same authors extended this result to the following theorem:
\begin{thm}\cite{bfmw}
For every $n\geq 6$ there exist non-resolvable generalized $n$-manifolds with the disjoint disks property.
\end{thm}
Theorem 3.3 is a corollary of the main result in \cite{bfmw} stating that every generalized $n$-manifold, $n\geq 6$, is the cell-like image of a generalized $n$-manifold possessing the disjoint disks property. This theorem was conjectured in \cite{bfmw2}.
Another conjecture from \cite{bfmw2} is the Homogeneity conjecture:
\begin{conj}\cite{bfmw2}
Every connected generalized $(n\geq 5)$-manifold satisfying the disjoint disks property is homogeneous.
\end{conj}
In Bryant \cite{br4} it was shown that a generalized $n$-manifold $X$, $n\geq 5$, having the disjoint disks property also satisfies a general position property for maps of polyhedra into $X$ that genuine $n$-manifolds possess. More precisely, if $P$ and $Q$ are polyhedra of dimensions $p$ and $q$, respectively, then all maps $f:P\to X$ and $g:Q\to X$ can be approximated by maps $f'$ and $g'$ such that
$(i)$ $\dim (f'(P)\cap g'(Q))\leq p+q-n$ and $(ii)$ $p+q-n\leq n-3$ implies $X\setminus (f'(P)\cap g'(Q))$ is $1-LCC$ in $X$. Moreover, if $2p+1\leq n$, then every map of $P$ into $X$ can be approximated by $1-LCC$ embeddings.
\begin{itemize}
\item Recall that a set $A\subset X$ is {\em $1$-locally co-connected in $X$} (br., $1-LCC$) if for every $x\in A$ and a neighborhood $U$ of $x$ in $X$ there is a neighborhood $V$ of $x$ in $X$ such that the inclusion induced homomorphism $\pi_1(V\setminus A)\to \pi_1(U\setminus A)$ is trivial. A map (or embedding) $f:Y\to X$ is said to be $1-LCC$ provided the set $f(Y)$ is $1-LCC$.
\end{itemize}

According to Theorem 3.3, a positive solution of Conjecture 3.4 would imply that the Bing-Borsuk conjecture is false for $n\geq 6$.
On the other direction, Bryant \cite{br1} suggested the following modification of the Bing-Borsuk conjecture:
\begin{conj}\cite{br1}$[$Modified Bing-Borsuk conjecture$]$
Every locally compact homogeneous $ANR$-space of dimension $n\geq 3$ is a generalized $n$-manifold.
\end{conj}
A partial result concerning the Modified Bing-Borsuk conjecture is an old result of Bredon \cite{bre}, reproved by Bryant \cite{br2}:
\begin{thm}\cite{bre},\cite{br2}
If $X$ is a locally compact homogeneous $ANR$-space of dimension $n$ such that the groups $H_k(X,X\setminus\{x\};\mathbb Z)$, $k\leq n$, are finitely generated, then $X$ is a generalized $n$-manifold.
\end{thm}
Another result related to the Modified Bing-Borsuk conjecture was obtained by Bryant \cite{br} answering a question of Quinn \cite{qu}:
\begin{thm}\cite{br}
Every $n$-dimensional homologically arc-homogeneous $ANR$-compactum is a generalized manifold.
\end{thm}

\begin{itemize}
\item Here, a space $X$ is homologically arc-homogeneous \cite{qu} if for every path $\alpha:\mathbb I=[0,1]\to X$ the inclusion induced map
$$H_*(X\times\{0\},X\times\{0\}-(\alpha(0),0))\to H_*(X\times\mathbb I,(X\times\mathbb I)-\Gamma(\alpha))$$
is an isomorphism, where $\Gamma(\alpha)$ is the graph of $\alpha$.
\end{itemize}

More information about generalized manifolds can be found in Bryant \cite{br3}.

The last two theorems in this section show that $\mathbb Z$-homology have also common properties with Euclidean manifolds.

\begin{thm}\cite{kr1}
Let $X$ be a locally compact, locally connected $\mathbb Z$-homology $n$-manifold with $\dim X=n>1$ at each point. Then $X$ is a local Cantor manifold, i.e. every open connected subset of $X$ is a Cantor manifold.
\end{thm}
This result was extended in \cite[Corollary 4.2]{tv}.
\begin{thm}\cite{tv} Let $X$ be a complete metric space which is a $\mathbb Z$-homology $n$-manifold. Then every open arcwise connected subset of $X$ is a Mazurkiewicz arc manifold with respect to the class of all spaces of dimension $\leq n-2$.
\end{thm}

\section{Local homological and cohomological structure of homogeneous $ANR$-compacta}

We are going to show in this section that the local cohomological and homological structure of homogeneous $n$-dimensional $ANR$-continua is similar to the corresponding local structure of $\mathbb R^n$.
\begin{itemize}
\item Recall that for any abelian group $G$ the cohomology group $\check{H}^n(X;G)$ is
isomorphic to the group of pointed homotopy classes of maps from $X$ into the Eilenberg-MacLane space $K(G,n)$ of type $(G,n)$, see cite{spa}. The cohomological dimension $\dim_GX$ is the largest integer $m$ such there exists a closed set $A\subset X$ with $\check{H}^m(X,A;G)\neq 0$. Equivalently, $\dim_GX\leq n$ if and only if every map $f:A\to K(G,n)$ can be extended to a map $\widetilde f:X\to K(G,n)$.
\end{itemize}

Suppose $(K,A)$ is a pair of closed subsets of a space $X$ with $A\subset K$. Then we denote by $j^n_{K,A}:\check{H}^n(K;G)\to\check{H}^n(A;G)$
and $i^n_{A,K}:\check{H}_n(A;G)\to\check{H}_n(K;G)$, respectively, the inclusion induced cohomology and homology homomorphisms. We say that
an element $\gamma\in\check{H}^n(A;G)$ is not extendable over $K$ if $\gamma$ is not contained in the image $j^n_{K,A}(\check{H}^n(A;G))$.
\begin{itemize}
\item If $(K,A)$ is as above, we say that $K$ is an {\em $(n,G)$-homology membrane spanned on $A$ for an element $\gamma\in\check{H}_n(A;G)$} provided $i^n_{A,K}(\gamma)=0$, but $i^n_{A,P}(\gamma)\neq 0$ for every proper closed subset $P$ of $K$ with $A\subset P$. Similarly,
$K$ is said to be an {\em $(n,G)$-cohomology membrane spanned on $A$ for an element $\gamma\in\check{H}^n(A;G)$} if $\gamma$ is not extendable over $K$, but it is extendable over any proper closed subset $P$ of $K$ containing $A$.
\end{itemize}
The continuity of the \v{C}ech cohomology \cite{es} implies the following fact: If $A$ is a closed subset of a compact space $X$ and $\gamma\in\check{H}_n(A;G)$ is not extendable over $X$, then there is an $n$-homology membrane for $\gamma$ spanned on $A$. We also note that under the same assumption for $A$ and $X$, the existence of a non-trivial $\gamma\in\check{H}_n(A;G)$ with $i^n_{A,X}(\gamma)=0$ yields the existence of a closed set $K\subset X$ containing $A$ such that $K$ is an $n$-homology membrane for $\gamma$ spanned on $A$, see \cite{bb}.

Next theorem provides the local cohomological structure of homogenous $ANR$s.
\begin{thm}\cite{vv1} Let $X$ be a homogeneous $ANR$-continuum with $\dim_GX=n\geq 2$ and $G$ be a countable principal ideal domain.
Then every point $x\in X$ has a basis $\mathcal B_x$ of open sets $U\subset X$ satisfying the following conditions:
\begin{itemize}
\item[(1)] $\rm{int}\overline U=U$ and the complement of $\bd\, U$ has exactly two components;
\item[(2)] $\check{H}^{n-1}(\bd\, U;G)\neq 0$, $\check{H}^{n-1}(\overline U;G)=0$ and $\overline{U}$ is an $(n-1,G)$-cohomology membrane spanned on $\bd\, U$ for any non-zero $\gamma\in\check{H}^{n-1}(\bd\, U;G)$;
\item[(3)] $\bd\, U$ is a cohomological $(n-1,G)$-bubble;
\item[(4)] The inclusion homomorphism $j_{U,V}^n:\check{H}^{n}(X,X\setminus U;G)\to \check{H}^{n}(X,X\setminus V;G)$ is non-trivial for any $U,V\in\mathcal B_x$ with $U\subset V$.
\end{itemize}
\end{thm}
\begin{itemize}
\item Here, a closed set $A\subset X$ is called a {cohomological \em $(n,G)$-bubble} if $\check{H}^n(A;G)\neq 0$ but $\check{H}^n(B;G)=0$ for every closed proper subset $B\subset A$. Considering in that definition \v{C}ech homology instead of cohomology, we obtain the notion of
    {homological \em $(n,G)$-bubble}.
\end{itemize}
\begin{cor}\cite{vv1} Let $X$ be a homogeneous $ANR$-compactum with $\dim_GX=n\geq 2$ and $G$ be a countable group. Then
\begin{itemize}
\item[(1)] $f(U)$ is open in $X$ provided $U\subset X$ is open and $f:U\to X$ is an injective map;
\item[(2)] $\dim_GA=n$, where $A\subset X$ is closed, if and only if $A$ has a non-empty interior in $X$.
\end{itemize}
\end{cor}
\begin{itemize}
\item For any abelian group $G$ Alexandroff \cite{ps1} introduced the dimension $d_GX$ of a space $X$ as the maximum integer $n$ such that there exists a closed set $F\subset X$ and a non-trivial $\gamma\in\check{H}_{n-1}(F;G)$ such that $i^{n-1}_{F,X}(\gamma)=0$. We have the following inequalities for any finite-dimensional metric compactum $X$ and any $G$: $d_GX\leq\dim X=d_{\mathbb Q_1}X=d_{\mathbb S^1}X$, where $\mathbb S^1$ is the circle group and $\mathbb Q_1$ is the group of rational elements of $\mathbb S^1$.
\end{itemize}
Because the definition of $d_GX$ does not provide any information for the homology groups $\check{H}_{k-1}(F;G)$ when $F\subset X$ is closed and $k<d_GX-1$, we consider
the set $\mathcal{H}_{X,G}$ of all integers $k\geq 1$ such that there exist a closed set $F\subset X$ and a non-trivial element $\gamma\in \check{H}_{k-1}(F;G)$ with $i^{k-1}_{F,X}(\gamma)=0$. Obviously, $d_GX=\max\mathcal{H}_{X,G}$.
\begin{thm}\cite{vv}
Let $X$ be a finite dimensional homogeneous metric $ANR$ with $\dim X\geq 2$. Then every point $x\in X$ has a basis $\mathcal B_x=\{U_k\}$ of open sets such that for any abelian group $G$ and $n\geq 2$ with $n\in\mathcal{H}_{X,G}$ and $n+1\not\in\mathcal{H}_{X,G}$ almost all $U_k$ satisfy the following conditions:
\begin{itemize}
\item[(1)] $\check{H}_{n-1}(\bd\overline U_k;G)\neq 0$ and $\overline{U}_k$ is an $(n-1,G)$-homology membrane spanned on $\bd\overline U_k$ for any non-zero $\gamma\in\check{H}_{n-1}(\bd\overline U_k;G)$;
\item[(2)] $\check{H}_{n-1}(\overline U_k;G)=\check{H}_{n}(\overline U_k;G)=0$ and $X\setminus\overline U_k$ is connected;
%\item[(3)] $H_{n-1}(bd\overline U;G)\neq 0$;
\item[(3)] $\bd\overline U_k$ is a homological $(n-1,G)$-bubble.
\end{itemize}
\end{thm}
\begin{itemize}
\item A closed set $F\subset X$ is called {\em strongly contractible in $X$} if there is a proper closed subset $A\subset X$ such that $F$ is contractible in $A$ and $A$ is contractible in $X$.
\end{itemize}

\begin{cor}\cite{vv}
Let $X$ be a homogeneous compact metrizable $ANR$-space such that $n\in\mathcal{H}_{X,G}$ and $n+1\not\in\mathcal{H}_{X,G}$. Then for every closed set $F\subset X$ we have:
\begin{itemize}
\item[(1)] $\check{H}_{n}(F;G)=0$ provided $F$ is contractible in $X$;
\item[(2)] $F$ separates $X$ provided $\check{H}_{n-1}(F;G)\neq 0$ and $F$ is strongly contractible in $X$;
\item[(3)] If $K$ is a homological membrane for some non-trivial element of $\check{H}_{n-1}(F;G)$ and $K$ is contractible in $X$, then $(K\setminus F)\cap\overline{X\setminus K}=\varnothing$.
\end{itemize}
\end{cor}

\section{Cyclicity, full-valuedness and existence of non-degenerate finite-dimensional homogeneous $AR$s}

In this section we discuss more problems about homogeneous $ANR$-compacta. We denote by $\mathcal H(n)$ the class of all homogeneous metric $ANR$-compacta of dimension $n$. If not explicitly stated otherwise, everywhere in this section $X$ is a space from  $\mathcal H(n)$.
\begin{qu}\cite{bb}$($Cyclicity$)$ Is it true that:
\begin{itemize}
\item[(1)] $X$ is cyclic in dimension $n$$?$
\item[(2)] No closed subset of $X$, acyclic in dimension $n-1$, separates $X$$?$
\end{itemize}
\end{qu}
\begin{itemize}
\item We say that $X$ is {\em cyclic in dimension $n$} if there is a group $G$ such that $\check{H}^n(X;G)\neq 0$. If a space is not cyclic in dimension $n$, it is called {\em acyclic in dimension $n$}. %Actually, the original cyclicity question $(1)$ was whether $\check{H}_n(X;\mathbb Z)\neq 0$.
    Because $X$ is $n$-dimensional $ANR$, the duality between \v{C}ech homology and cohomology (see \cite{hw}) and the universal coefficient formulas imply the following equivalence: $\check{H}_n(X;G)\neq 0$ for some group $G$ iff and only if $X$ cyclic in dimension $n$.
\end{itemize}
 Next results show that the two parts of  the cyclicity Question 5.1 have simultaneously positive or negative answers.
\begin{thm}\cite{vv}
The following conditions are equivalent:
\begin{itemize}
\item[(1)] For all $n\geq 1$ and $X\in\mathcal H(n)$ there exists a group $G$ with $\check{H}^n(X;G)\neq 0$ $($resp., $\check{H}_n(X;G)\neq 0)$;
\item[(2)] If $X\in\mathcal H(n)$, $n\geq 1$, and $F\subset X$ is a closed set separating $X$, then there exists a group $G$ with $\check{H}^{n-1}(F;G)\neq 0$ $($resp., $\check{H}_{n-1}(F;G)\neq 0)$;
\item[(3 )] If $X\in\mathcal H(n)$, $n\geq 1$, and $F\subset X$ is a closed set separating $X$ with $\dim F\leq n-1$, then there exists a group $G$ such that $\check{H}^{n-1}(F;G)\neq 0$ $($resp., $\check{H}_{n-1}(F;G)\neq 0)$.
\end{itemize}
\end{thm}

On the other hand, the structure of cyclic homogeneous $ANR$ continua is described in \cite{vv3} (see sections 2 and 4, respectively, for the definitions of a strong $V^n_G$-continuum and a cohomological $(n,G)$-bubble).
\begin{thm}\cite{vv3}
Let $X$ be a homogeneous metric $ANR$-continuum such that $\dim_GX=n$ and $\check{H}^n(X;G)\neq 0$ for some group $G$. Then
\begin{itemize}
\item[(1)] $X$ is a cohomological $(n,G)$-bubble;
\item[(2)] $X$ ix a strong $V^n_G$-continuum;
\item[(3)] $\check{H}^{n-1}(A;G)\neq 0$ for every closed set $A\subset X$ separating $X$.
\end{itemize}
\end{thm}
Items $(1)$ and $(3)$ were also established by Yokoi \cite{yo} for the case $G$ is a principal ideal domain.

Actually, the third item of Theorem 5.3 is true for every set $A\subset X$ cutting $X$ between two disjoint open subsets of $X$, see \cite{tv}.
Recall that {\em $A$ cuts $X$ between two disjoint sets $U$ and $V$} if $A\cap (U\cup V)=\varnothing$ and every continuum in $X$ joining $U$ and $V$ meets $A$.

Bing-Borsuk \cite{bb} proved that all locally homogeneous and locally compact $ANR$ spaces of dimension $n=0,1,2$ are $n$-manifolds. Therefore,
there is no such a non-degenerated $AR$-space of dimension $\leq 2$. So, next question is interesting only in dimension $\geq 3$.
\begin{qu}\cite{bb}, \cite{bo} Does there exists non-degenerate finite-dimensional homogeneous $AR$-compactum$?$
\end{qu}
%Both of the above two questions were formulated for locally homogeneous spaces.
Another questions in that direction was listed in\cite{west}.
\begin{qu}\cite{west} Is the Hilbert cube $Q$ the only homogeneous non-degenerate compact $AR$$?$
\end{qu}
Obviously, if every finite-dimensional homogeneous $ANR$ is cyclic, then the answer of Question 5.4 is no.
Assuming that there exists a non-degenerate homogeneous finite-dimensional $AR$-space $X$, a possible way to obtain a contradiction is to find a continuous map $f:X\to X$ without having fix points. This is equivalent to find a continuous selection of the set-valued map $\varphi:X\rightsquigarrow X$, $\varphi(x)=X\setminus\{x\}$. To find such a selection we need to know that all sets $X\setminus\{x\}$ have
nice properties of type $C^{n-1}$, and that is not clear. More refined homological selection theorems (see \cite{bc} or \cite{vv2}) could be helpful. For example, in case $G$ is a field, the following result could be useful:
\begin{thm}\cite{vv2}
Let $X$ be a compact metric $AR$ and $\Phi:X\rightsquigarrow X$ be an upper semi-continuous compact-valued homological $UV^{n-1}(G)$-map. Then $\Phi$ has a fixed point, i.e. there is $x_0\in X$ with $x_0\in\Phi(x_0)$.
\end{thm}
\begin{itemize}
\item Here, a set a closed set $A\subset X$ is said to be {\em homological $UV^k(G)$} \cite{vv2} if for every neighborhood $U$ of $A$ in $X$ there exists another neighborhood $V$ of $A$ such that $V\subset U$ and all inclusion homomorphisms $\check{H}_m(V;G)\to \check{H}_m(U;G)$, $m\leq k$,  are trivial. A compact-valued map $\Phi:X\rightsquigarrow X$ is a {\em homological $UV^k(G)$-map} provided each $\Phi(x)$ is a
    homological $UV^k(G)$-set in $X$.
\end{itemize}

Next question goes back to Bryant \cite{br1}, it was also discussed in \cite{clqr} and \cite{f}.
\begin{qu}
Is it true that any homogeneous finite-dimensional $ANR$-compactum is dimensionally full-valued$?$
\end{qu}
\begin{itemize}
\item A compactum $X$ is {\em dimensionally full-valued} if $\dim (X\times Y)=\dim X+\dim Y$ for every compact space $Y$.
\end{itemize}
Dranishnikov \cite{dra} constructed a family of $4$-dimensional metric $AR$-compacta $M_p$, where $p$ is a prime number, such that $\dim (M_p\times M_q)=7$ for all $p\neq q$. The spaces $M_p$ are not homogeneous. On the other hand, the classical Pontryagin surfaces \cite{po} is a family of $2$-dimensional homogeneous metric, but not $ANR$-compacta $\{\Pi_p:p{~}\mbox{is prime}\}$  with $\dim (\Pi_p\times\Pi_q)=3$ for $p\neq q$.

Boltyanskii \cite{bol} provided a criterion for dimensional full-valuedness:
\begin{thm}\cite{bol}
A finite-dimensional compactum is dimensionally full-valued if $\dim_GX=\dim X$ for any group $G$.
\end{thm}
Because $\dim_{\mathbb Q}X\leq\dim_GX$ for any $ANR$-compactum $X$ and all groups $G$ (for example, see \cite{dra1}), where $\mathbb Q$ is the group of rationals, it follows that a finite-dimensional $ANR$-compactum is dimensionally full-valued iff $\dim_{\mathbb Q}X=\dim X$.
Another criterion for $ANR$-space is established in \cite{vv1}:
\begin{thm}\cite{vv1}
The following conditions are equivalent for every $X\in\mathcal H(n)$:
\begin{itemize}
\item[(1)] $X$ is dimensionally full-valued;
\item[(2)] There exists a point $x\in X$ with $\check{H}_n(X,X\setminus\{x\};\mathbb Z)\neq 0$;
\item[(3)] $\dim_{\mathbb S^1}X=n$.
\end{itemize}
\end{thm}
Since $\check{H}_3(x,X\setminus\{x\};\mathbb Z)\neq 0$ for every $x\in X$ \cite{kr1}, where $X$ is a homogeneous metric $ANR$-compactum with $dim X=3$, Theorem 5.9 implies the following:
\begin{cor}\cite{vv1}
Every homogeneous metric $ANR$-compactum of dimension $3$ is dimensionally full-valued.
\end{cor}
Theorem 5.9 also implies that every compact generalized $n$-manifold is dimensionally full-valued. More generally, a positive answer of the following question yields that every $X\in\mathcal H(n)$ is dimensionally full-valued:
\begin{qu}[A weaker version of the modified Bing-Borsuk conjecture]
Is it true that $\check{H}_k(X,X\setminus\{x\};\mathbb Z)=0$ for every $x\in X$ and $k\leq n-1$, where $X\in\mathcal H(n)$$?$
\end{qu}

We still don't know if the dimension of any product of two homogeneous $ANR$-compacta obeys the logarithmic law.
\begin{qu}\cite{clqr}
Does the equality $\dim (X\times Y)=\dim X+\dim Y$ hold for any homogeneous $ANR$-compacta $X$ and $Y$$?$
\end{qu}
%%%%%%%%%% Bibliography %%%%%%%%%%%%%%%%%%%%%%%%%%


\begin{thebibliography}{999}

\bibitem{ps1} P.~Alexandroff. Introduction to homological dimension theory and general combinatorial topology.
Moscow, Nauka, 1975 (in Russian).

\bibitem{ps} P.~Alexandroff. Die Kontinua $(V^p)$ - eine
Versch\"{a}rfung der Cantorschen Mannigfaltigkeiten. \textit{Monatshefte
fur Math.} \textbf{61} (1957), 67--76 (German).

%\bibitem{bck}
%T.~Banakh,~R.~Cauty and A.~Karassev, \textit{On homotopical homological $Z_n$-sets}, Top. Proc. \textit{38} (2011), 29--82.

\bibitem{bc}
T.~Banakh and R.~Cauty. A homological selection theorem implying a division theorem for $Q$-manifolds.
\textit{Fixed point theory and its applications}, 11--22, Banach center Publ. \textit{77}, Polish Acad. Sci. Inst. Math., Warsaw, 2007.

\bibitem{bdvw}
M.~Bestvina, R.~Daverman, G.~Venna and J.~Walsh. A $4$-dimensional $1-LCC$ shrinking theorem. \textit{Topology Appl.} \textbf{110} (2001), 3--20.

%\bibitem{bi}
%R.~H.~Bing, \textit{Partitioning a set}, Bull. Amer. Math. Soc. \textbf{55} (1949), no. 12, 1101--1110.

\bibitem{bb}
R.~H.~Bing and K.~Borsuk. Some remarks concerning topological homogeneous spaces.  \textit{Ann. of Math.} \textbf{81}, 1 (1965), 100--111.

\bibitem{bo}
K.~Borsuk. Theory of retracts. \textit{Monografie Matematyczne} \textbf{44}, PWN, Warsaw, 1967.

\bibitem{bol}
V.~Boltyanskii. On dimensional full-valuedness of compacta. \textit{Dokl. Akad. Nauk SSSR} \textbf{67} (1949), 773--777 (in Russian).

\bibitem{bre}
G.~Bredon. Sheaf Theory. Second Edition, Graduate texts in Mathematics 170, Springer, 1997.

\bibitem{br}
J.~Bryant. Homologically arc-homogeneous $ENR$s. \textit{Geometry and Topology Monographs} \textbf{9} (2006), 1--6.

\bibitem{br1}
J.~Bryant. Reflections on the Bing-Borsuk conjecture. Abstracts of the talks presented at the 19th Annual Workshop in Geometric Topology,
June 13-15, 2002, 2-3.

\bibitem{br3}
J.~Bryant. A survey of recent results on generalized manifolds. \textit{Topology Apll.} \textbf{113} (2001), 13--22.

\bibitem{br2}
J.~Bryant. Homogeneous $ENR$'s. \textit{Topology Appl.} \textbf{27} (1987), 301--306.

\bibitem{br4}
J.~Bryant. General position theorems for generalized manifolds. \textit{Proc. Amer.Math. Soc.} \textbf{98}, 4 (1986), 667--670.

\bibitem{bfmw}
J.~Bryant, S.~Ferry, W.~Mio and S.~Weinberger. Desingularizing homology manifolds. \textit{Geom. Topol.} \textbf{11} (2007), 1289--1314.

\bibitem{bfmw1}
J.~Bryant, S.~Ferry, W.~Mio and S.~Weinberger. Topology of homology manifolds.  \textit{Ann. Math.} \textbf{143} (1996), 435--467.

\bibitem{bfmw2}
J.~Bryant, S.~Ferry, W.~Mio and S.~Weinberger. Topology of homology manifolds. \textit{Bull. Amer. Math. Soc. (N.S.)} \textbf{28}, 2 (1993),  324–328.

\bibitem{bro}
L.~Brouwer. Zur Invarianz des $n$-dimensionalen Gebiets. \textit{Math. Annalen} \textbf{72} (1912), 55--56.

\bibitem{bu1}
H.~Busemann. Metric methods in Finsler spaces and in the foundation of geometry. \textit{Ann. Math. Study} \textbf{8}, Princeton University Press,
Princeton, 1942.

\bibitem{bu2}
H.~Busemann. On spaces in which two points determine a geodesic. \textit{Trans. Amer. Math. Soc.} \textbf{54} (1943), 171--184.

\bibitem{can1}
J.~Cannon. Shrinking cell-like decompositions of manifolds. Codimension three.  \textit{Ann. of Math.} \textbf{110}, 1 (1979), 83–112.

\bibitem{can}
J.~Cannon. The recognition problem: what is a topological manifold. \textit{Bull. Amer. Math. Soc.} \textbf{84} (1978), 832--866.

\bibitem{clqr}
M.~Cardenas, F.~Lasheras,A.~Quintero and D.~Repov\v{s}. On manifolds with nonhomogeneous factors. \textit{Centr. Eur. J. Math.} \textbf{10}
(2012), 857--862.

%\bibitem{ca}
%R.~Cauty, \textit{La classe bor\'{e}lienne ne d\'{e}termine pas le type topologique de $C_p(X)$}, Serdica Math. J. \textbf{24} (1998), 307--318.

%\bibitem{ch}
%J.~Choi, \textit{Properties of $n$-bubles in $n$-dimensional compacta and the existence of $(n-1)$-bubles in
%$n$-dimensional $clc^n$ compacta}, Top. Proceed. \textbf{23} (1998), 101-120.

\bibitem{dr1}
R.~Daverman and D.~Repov\v{s}. General position properties that characterize $3$-manifolds. \textit{Canad. J. Math.} \textbf{44} (1992), 234--251.

\bibitem{dr}
R.~Daverman and D.~Repov\v{s}. A new $3$-dimensional shrinking criterion. \textit{Trans. Amer. Math. Soc.} \textbf{315} (1989), 219--230.

\bibitem{dra1},
A.~Dranishnikov. Cohomological dimension theory of compact metric spaces. Topology Atlas invited contribution, vol. \textbf{6}, 2001, 7--73

\bibitem{dra}
A.~Dranishnikov. Homological dimension theory. \textit{Russian Math. Surveys} \textbf{43}, 4 (1988), 11--63.
%\bibitem{dr} T.~Dobrowolski and L.~Rubin, \textit{The hyperspaces of
%infinite-dimensional compacta for covering and cohomological
%dimension are homeomorphic}, Pacif. J.Math. \textbf{164} (1994), 15--39.

%\bibitem{dm}
%J.~Dugundji and E.~Michael, \textit{On local and uniformly local topological properties}, Proc. Amer. Math. Soc. \textbf{7} (1956), 304--307.

%\bibitem{d}
%J.~Dugundji, \textit{Modified Vietoris theorems for homotopy}, Fund. Math. \textbf{66} (1969), 223--235.

\bibitem{ed}
R.~Edwards. The topology of manifolds and cell-like maps. \textit{Proceedings of the International Congress of Mathematicians} (Helsinki, 1978), Acad. Sci. Fennica, Helsinki, 1980, 111-127.

%\bibitem{ef} E. G. Effros, \textit{Transformation groups and $C^*$-algebras},
%Ann. of Math. \textbf{81} (1965), 38--55.

\bibitem{es}
S.~Eilenberg, N.~Steenrod. Foundations of algebraic topology. Proceton Univ. Press, Prinstone, New Jersey, 1952.

\bibitem{en} R.~Engelking. Theory of
dimensions: Finite and Infinite. Heldermann Verlag, Lemgo, 1995.

\bibitem{f}
V.~Fedorchuk. On homogeneous Pontryagin surfaces. \textit{Dokl. Akad. Nauk} \textbf{404}, 5 (2005),  601--603 (in Russian).

\bibitem{ja}
W.~Jaconsche. The Bing-Borsuk conjecture is stronger than the Poincare conjecture. \textit{Fund. Math.} \textbf{106} (1980), 127--134.

\bibitem{ht}
N.~Hadjiivanov, V.~Todorov. On non-Euclidean manifolds. \textit{C. R. Acad. Bulgare Sci.} \textbf{33} (1980), 449--452 (in Russian).

\bibitem{hr}
D.~Halverson and D.~Repov\v{s}. The Bing-Borsuk conjecture and the Busemann conjecture. \textit{Math. Communications} \textbf{13} (2008), 163--184.

%\bibitem{re:95} R.~Engelking, \textit{Theory of
%dimensions: Finite and Infinite}, Heldermann Verlag, Lemgo, 1995.

\bibitem{hw} W. Hurewicz and H. Wallman. Dimension theory.
Princeton University Press, Princeton, 1948.

%\bibitem{hu}
%P.~Huber, \textit{Homotopical cohomology and \v{C}ech cohomology}, Math. Annalen \textbf{144} (1961), 73--76.

\bibitem{ka}
A.~Karassev, A private communication.

\bibitem{ktv}
A.~Karassev, V.~Todorov and V.~Valov. Alexandroff manifolds and homogeneous continua. \textit{Canad. Math. Bull.} \textbf{57}, 2 (2014), 335--343.

\bibitem{kktv}
A.~Karassev, P.~Krupski, V.~Todorov and V.~Valov. Generalized Cantor manifolds and homogeneity. \textit{Housto J. Math.} \textbf{38}, 2 (2012), 583--609.

\bibitem{kr1}
P.~Krupski. On the disjoint $(0,n)$-cells property for homogeneous ANRs. \textit{Colloq. Math.} \textbf{66}, 1 (1993), 77—84.

%\bibitem{kr2}
%P.~Krupski, \textit{Recent results on homogeneous curves and ANRs}, Topology Proc. \textbf{16} (1991), 109—118.

\bibitem{kr3}
P.~Krupski. Homogeneity and Cantor manifolds. \textit{Proc. Amer. Math. Soc.} \textbf{109} (1990), 1135--1142.

\bibitem{kv}
P.~Krupski and V.~Valov. Mazurkiewicz manifolds and homogeneity. \textit{Rocky J. Math.} \textbf{41}, 6 (2011), 1933--1938.

%\bibitem{ku} V.~Kuz'minov, \textit{Homological dimension theory}, Russian Math. Surveys \textbf{23} (1968), no. 1, 1--45.

\bibitem{ly}
J.~Lysko. Some theorems concerning finite dimensional homogeneous ANR-spaces. \textit{Bull. Acad. Polon. Sci. Sér. Sci. Math. Astronom. Phys.}
\textbf{24}, 7 (1976), 491–496.

\bibitem{per}
J.~Morgan and G.~Tian. Ricci flow and the Poincare conjecture. \textit{Clay Math. Monographs} \textbf{3}, Amer. Math. Soc., Providence, RI; Clay Math. Institute, Cambridge, MA, 2007.

\bibitem{po}
L.~Pontryagin. Sur une hypothese foundamentale de la dimension. \textit{C.R. Acad. Sci.} \textbf{190} (1930), 1105--1107.

\bibitem{qu}
F.~Quinn. Problems on homology manifolds. \textit{Geometry and Topology Monographs} \textbf{9} (2006), 87--103.

\bibitem{re}
D.~Repov\v{s}. The recognition problem for topological manifolds. in: \textit{Geometric and Algebraic Topology}, (J. Krasinkiewicz, S. Spiez, and H. Torunczyk, Eds.), PWN, Warsaw, 1986, 77--108.

\bibitem{se}
H.~Seidel. Locally homogeneous $ANR$-spaces. \textit{Arch. Math.} \textbf{44} (1985), 79--81.

%\bibitem{spa}
%E.~Spanier, \textit{Algebraic Topology}, McGraw-Hill Book Company, 1966.

\bibitem{s}
S.~Stefanov. A cohomological analogue of $V^n$-continus and a theorem of Mazurkiewicz. \textit{Serdica Math. J.} \textbf{12}, 1 (1986), 88--94.

\bibitem{t}
V.~Todorov. Irreducibly cyclic compacta are Cantor manifolds. \textit{Proc. of the Tenth Spring Conf. of the Union of Bulg. Math.} Sunny Beach, April 6-9 1981.

\bibitem{tv}
V.~Todorov and V.~Valov. Alexandroff type manifolds and homology manifolds. \textit{Houston J. Math.} \textbf{40}, 4 (2014), 1325--1346.

\bibitem{ur}
P.~Urysohn. Memoire sur les multiplicites cantoriennes. \textit{Fund. Math.} \textbf{7} (1925), 30--137.
%\bibitem{tv}
%V.~Todorov and V.~Valov, \textit{Generalized Cantor manifolds and indecomposable continua},
%Questions and Answers in Gen. Topolol., accepted.

\bibitem{yo}
K.~Yokoi. Bubbly continua and homogeneity. \textit{Houston J. Math.} \textbf{29}, 2 (2003), 337--343.

\bibitem{vv}
V.~Valov. Local homological properties and cyclicity of homogeneous $ANR$ compacta.
\textit{Proc. Amer. Math. Soc.} \textbf{146} (2018), 2697--2705.

\bibitem{vv2}
V.~Valov. Homological selections and fixed-point theorems. \textit{J. Fixed Point Theory} \textbf{19} (2017), 1561--1570.

\bibitem{vv1}
V.~Valov. Local cohomological properties of homogeneous $ANR$ compacta. \textit{Fund. Math.} \textbf{223} (2016), 257--270.

\bibitem{vv3}
V.~Valov. Homogeneous $ANR$-spaces and Alexandroff manifolds. \textit{Topology Appl.}
\textbf{173} (2014), 227--233.

\bibitem{west}
J.~West. Open problems in infinite-dimensional topology. in: \textit{Open Problems in Topology} (J. van Mill and G.M. Reed, Eds.),
Elsevier Science Publishers B. V., North-Holland, 1990, 524--597.

\bibitem{w}
R.~Wilder. Topology of manifolds. \textit{Amer. Math. Soc. Coll.} \textbf{32} (1949).
\end{thebibliography}
\end{document}